\documentclass[lettersize,journal]{IEEEtran}
\usepackage{amsmath,amsfonts}
\usepackage{algorithmic}
\usepackage{algorithm}
\usepackage{array}
\usepackage[caption=false,font=normalsize,labelfont=sf,textfont=sf]{subfig}
\usepackage{textcomp}
\usepackage{stfloats}
\usepackage{multicol}
\usepackage{adjustbox}
\usepackage[symbol]{footmisc}
\usepackage{booktabs} 
\usepackage{nomencl}
\usepackage{url}
\usepackage{verbatim}
\usepackage{graphicx}
\usepackage{cite}
\usepackage{siunitx}
\begin{document}

\title{Applications of Lifted Nonlinear Cuts to Convex Relaxations of the AC Power Flow Equations}

\author{Sergio I. Bugosen, Robert B. Parker, Carleton Coffrin\\
Los Alamos National Laboratory, Los Alamos, NM, USA}

\markboth{}
{Bugosen {\it et al.}: Applications of Lifted Nonlinear Cuts}


\maketitle

\begin{abstract}
We demonstrate that valid inequalities, or lifted nonlinear cuts (LNC), can be projected
to tighten the Second Order Cone (SOC), Convex DistFlow (CDF), and Network Flow (NF)
relaxations of the AC Optimal Power Flow (AC-OPF) problem. We conduct experiments on
38 cases from the PGLib-OPF library,
showing that the LNC strengthen the SOC and CDF relaxations in 100\% of the test cases,
with average and maximum differences in the optimality gaps of 6.2\% and 17.5\% respectively.
The NF relaxation is strengthened in 46.2\% of test cases, with average and maximum
differences in the optimality gaps of 1.3\% and 17.3\% respectively. We also study
the trade-off between relaxation quality and solve time, demonstrating that the strengthened
CDF relaxation outperforms the strengthened SOC formulation in terms of runtime and
number of iterations needed, while the strengthened NF formulation is the most scalable
with the lowest relaxation quality improvement due to these LNC.
\end{abstract}


\makenomenclature
\renewcommand\nomlabelwidth{2.5cm}
\nomenclature[\_1]{$(\cdot)^*$}{Conjugate of a complex number}
\nomenclature[\_2]{$\vert\cdot\vert$}{Euclidean magnitude of a complex number}
\nomenclature[\_3]{$\Re(\cdot)$}{Real part of a complex number}
\nomenclature[\_4]{$\Im(\cdot)$}{Imaginary part of a complex number}
\nomenclature[0]{$N$}{The set of nodes in the network}
\nomenclature[1]{$E$}{The set of \textbf{from} edges in the network}
\nomenclature[2]{$I$}{AC current}
\nomenclature[4]{$L$}{Current magnitude squared, $\vert I\vert^2$}
\nomenclature[3]{$V = v\angle{\theta}$}{AC voltage}
\nomenclature[5]{$S = p + \textbf{i}q$}{AC power}
\nomenclature[6]{$Y = g + \textbf{i}b$}{Line admittance}
\nomenclature[7]{$Z = r + \textbf{i}x$}{Line impedance}
\nomenclature[8]{$W_{ij} = w_{ij}^R + \textbf{i}w_{ij}^I$}{Product of two AC voltages $V_i$ and $V_j^*$}
\nomenclature[81]{$W_i$}{Voltage magnitude squared}
\nomenclature[9]{$\theta_{ij}$}{Phase angle difference (i.e., $\theta_i - \theta_j)$}
\nomenclature[90]{$\phi_{ij}, \delta_{ij}$}{Phase angle difference center and offset}
\nomenclature[941]{$x^l,x^u$}{Lower, upper bounds of a variable $x$}
\nomenclature[95]{$x^\sigma$}{Sum of the bounds (i.e., $x^l + x^u$)}
\nomenclature[96]{$\textbf{x}$}{A constant value}

\printnomenclature

\section{Introduction}
\IEEEPARstart{T}{he} AC Optimal Power Flow problem (AC-OPF) is fundamental in power
systems computations. It seeks to determine the operating conditions of an electric
network such that an objective function (often generation cost minimization) is optimized,
electricity demand is met, and AC power flow equalities are satisfied. This problem
contains nonconvex and nonlinear constraints, and is known to be NP-hard \cite{Lehmann2016-hn}.

Convex relaxations such as the Semi-definite Programming (SDP), Second Order Cone (SOC),
Convex DistFlow (CDF), Quadratic Convex (QC) and Network Flow (NF) formulations are
useful to provide bounds on the AC-OPF objective function, prove infeasibility of particular
instances, and produce a solution that, if found feasible in the original nonconvex
problem, guarantees that it is a global optimum \cite{Coffrin2017-fx}. Convex relaxations
are also useful to provide bounds in contexts where using a nonconvex model is intractable. 
Strengthened convex relaxations provide better performance in global optimization algorithms
by reducing the number of partitions required in branch-and-bound, or reducing the
number of iterations needed in multi-tree methods \cite{Chen2016, Liu2019, Nagarajan2019-pd}.



Convex relaxations must balance solution quality (tightness) with tractability. Coffrin
et al. \cite{Coffrin2017-fx,coffrin2016strengthening} develops a novel approach to derive lifted nonlinear cuts
for the AC power flow equations, specifically to strengthen the SDP and QC relaxations,
without significantly increasing solve time. In this paper, we extend the lifted nonlinear
cuts to the SOC \cite{Jabr2006}, CDF \cite{Farivar2011} and NF \cite{Coffrin2016} relaxations.
We demonstrate the improved quality of the relaxations and show the trade-off between
relaxation quality and solve time that exists among the tightened versions of these
three formulations. The computational study is conducted on 38 test cases from the
PGLib-OPF benchmark library \cite{babaeinejadsarookolaee2021power}, which features
realistic datasets incorporating bus shunts, line charging, and transformers.




\section{Strengthening Convex Relaxations}

The AC-OPF problem is NP-hard due to the nonconvex product of voltage variables $V_iV_j^*$.
This product can be lifted into a higher-dimensional space (i.e. the $W$-space), where
voltage phase information is lost. The absolute square of the voltage product is then
relaxed (Eq. (\ref{eq0}d)) to obtain the basis for the SOC, CDF, and NF relaxations, 

\begin{subequations}
\label{eq0}
\begin{equation}
    w_i = |V_i|^2~~\forall i \in N
\end{equation}
\begin{equation}
    W_{ij} = V_iV_j^*~~\forall (i, j) \in E
\end{equation}
\begin{equation}
    |W_{ij}|^2 = w_iw_j~~\forall (i, j) \in E
\end{equation}
\begin{equation}
    |W_{ij}|^2 \leq w_iw_j~~\forall (i, j) \in E
\end{equation}
\end{subequations}


Coffrin et al. \cite{Coffrin2017-fx,coffrin2016strengthening} propose a novel
approach to derive valid inequalities
in the $W$-space. These valid inequalities, referred as lifted nonlinear cuts (LNC),
have been proven to strengthen the SDP and QC relaxations. The LNC are shown in Eqs.
(\ref{eq1})-(\ref{eq2}), where
$\boldsymbol{\phi_{ij}}=(\boldsymbol{\theta_{ij}^u} + \boldsymbol{\theta_{ij}^l})/2$
and $\boldsymbol{\delta_{ij}}=(\boldsymbol{\theta_{ij}^u} - \boldsymbol{\theta_{ij}^l})/2$. 

\begin{equation}
\label{eq1}
\begin{split}
    &\boldsymbol{v_i^\sigma} \boldsymbol{v_j^\sigma} (w_{ij}^R\cos{\boldsymbol{\phi_{ij}}}+w_{ij}^I\sin{\boldsymbol{\phi_{ij}}}) \\
    &- \boldsymbol{v_j^u} \cos(\boldsymbol{\delta_{ij}}) \boldsymbol{v_j^\sigma} w_i - \boldsymbol{v_i^u} \cos(\boldsymbol{\delta_{ij}}) \boldsymbol{v_i^\sigma} \dfrac{(w_{ij}^R)^2 + (w_{ij}^I)^2}{w_i}  \\
    &\geq \boldsymbol{v_i^u} \boldsymbol{v_j^u} \cos(\boldsymbol{\delta_{ij}}) \times (\boldsymbol{v_i^l v_j^l} - \boldsymbol{v_i^u v_j^u})~~\forall (i, j) \in E
\end{split}
\end{equation}

\begin{equation}
\label{eq2}
\begin{split}
    &\boldsymbol{v_i^\sigma v_j^\sigma} (w_{ij}^R\cos{\boldsymbol{\phi_{ij}}}+w_{ij}^I\sin{\boldsymbol{\phi_{ij}}}) \\
    &- \boldsymbol{v_j^l} \cos(\boldsymbol{\delta_{ij}}) \boldsymbol{v_j^\sigma} w_i - \boldsymbol{v_i^l} \cos(\boldsymbol{\delta_{ij}}) \boldsymbol{v_i^\sigma} \dfrac{(w_{ij}^R)^2 + (w_{ij}^I)^2}{w_i} \\
    &\geq \boldsymbol{v_i^l v_j^l} \cos(\boldsymbol{\delta_{ij}}) \times (\boldsymbol{v_i^u v_j^u} - \boldsymbol{v_i^l v_j^l})~~ \forall (i, j) \in E
\end{split}
\end{equation}


These LNC are nonlinear, but can be linearized by lifting them to the $\mathbb{R}^4$
space \{$w_i$, $w_j$, $w_{ij}^R$, $w_{ij}^I$\} using Eq. (\ref{eq0}c). The goal of
this work is to project these LNC into the variable space of the CDF and NF relaxations,
and demonstrate that they provide tighter optimality gaps. Note that the LNC are by
default expressed in the $W$-space, thus they are directly applicable to strengthen
the SOC relaxation. 
To highlight the effectiveness of these LNC, we run an optimization-based bound tightening
(OBBT) algorithm for the voltage ($v_i$) and phase angle difference ($\theta_{ij}$)
variables using the QC relaxation \cite{Coffrin2015, qcrel}. The LNC benefits from
these procedure as they are derived using the bounds on these variables.


\subsection{Strengthened NF relaxation}

The voltage product defined as $W_{ij} = w_{ij}^R + \boldsymbol{i}w_{ij}^I$ is not
a variable in the NF relaxation. Instead, this formulation is defined in the space
of the following variables: \{$W_i$, $S_{ij}$\}. The AC line flow equation, solved
for the voltage product term, yields $W_{ij} = w_i - \boldsymbol{Z_{ij}}^*S_{ij}$; this equation
is the basis to derive expressions for $w_{ij}^R$ and $w_{ij}^I$ in terms of the NF
variables. These are shown in Eqs. (\ref{eq3})-(\ref{eq4}), and are used to replace
$w_{ij}^R$ and $w_{ij}^I$ in Eqs. (\ref{eq1})-(\ref{eq2}).

\begin{equation}
\label{eq3}
w_{ij}^R = \Re(w_i - \boldsymbol{Z_{ij}^*}S_{ij}) ~~ \forall (i, j) \in E
\end{equation}

\begin{equation}
\label{eq4}
w_{ij}^I = \Im(w_i - \boldsymbol{Z_{ij}^*}S_{ij}) ~~ \forall (i, j) \in E
\end{equation}

\subsection{Strengthened CDF relaxation}

This relaxation is defined in the space of the following variables: \{$W_i$, $L_{ij}$,
$S_{ij}$\}. The expression for $w_{ij}^R$ in terms of the CDF variables is shown in
Eq. (\ref{eq5}) and is obtained by computing the absolute square of the AC current,
namely $L_{ij} = I_{ij}I_{ij}^* = |\boldsymbol{Y_{ij}}|^2(w_i - W_{ij} - W_{ij}^* + w_j)$. The expression
for $w_{ij}^I$ is equivalent to Eq. (\ref{eq4}). These equations are meant to replace
$w_{ij}^R$ and $w_{ij}^I$ in Eqs. (\ref{eq1})-(\ref{eq2}). Even though Eq. (\ref{eq3})
is also in the variable space of the CDF relaxation, preliminary experiments demonstrated
that the inclusion of the $L_{ij}$ variable in $w_{ij}^R$ is necessary to improve the
runtime performance of this formulation. 

\begin{equation}
\label{eq5}
w_{ij}^R = \frac{1}{2}\bigg(w_i + w_j - \frac{L_{ij}}{\vert \boldsymbol{Y_{ij}}\vert^2}\bigg)~~\forall (i, j) \in E
\end{equation}

\section{Computational Evaluation}

This section presents the benefits of strengthening the SOC, CDF, and NF relaxations
with their associated LNC projections, which were extended and implemented with bus
shunts, line charging, and transformers. The formulations for the SOC and CDF relaxations
can be found in \cite{coffrin2018distflow}, while the formulation for NF is in \cite{Coffrin2016}.
Since the LNC are an upper bound on branch line losses, we present results for the objective
of maximizing real power generation. These types of problems are present in a range of applications,
such as robust optimization \cite{Molzahn2018} and determination of voltage stability
margins \cite{Molzahn2013}.

Table \ref{table1} presents a comparison of optimality gaps and solve times for small instances from
PGLib-OPF, with and without the LNC. These instances have been preprocessed using OBBT.
Table \ref{table2} presents optimality gaps and solve times for large instances from PGLib-OPF
with the LNC applied. Due to long runtimes required to perform OBBT on large data sets, these
instances were not preprocessed with OBBT.
In Tables \ref{table1} and \ref{table2}, the SOC results were obtained with IPOPT v3.14
\cite{ipopt}, while the CDF and NF results were obtained with Gurobi v11.0 \cite{gurobi}.
These are the solvers that solved fastest, on average, for the respective formulations.
Optimality gaps are computed using a locally optimal AC-feasible solution as a lower bound
for the solution to the maximization problem. For instances in Table \ref{table1}, local
solutions were computed using IPOPT. For instances in Table \ref{table2}, local solutions
were computed using Knitro v14.0 \cite{knitro} as it was found to converge faster for
large instances.
All models were constructed using JuMP v1.23 \cite{Lubin2023} and PowerModels v0.21
\cite{powermod}.
Computations were conducted on a machine with an Apple M1 Max processor and 32 GB of
RAM running MacOS v13.6. The results presented in Tables \ref{table1} and \ref{table2}
are summarized below.

\begin{table*}[t]
\captionsetup{justification=centering} 
    \caption{Optimality gaps and runtime results for the power generation maximization
    problem. Test cases preprocessed with OBBT.}
\centering
\begin{adjustbox}{max width=\textwidth}

\begin{tabular}{|l|*{6}{S}|*{6}{S}|} \hline 
& \multicolumn{6}{c|}{\% Optimality Gap} & \multicolumn{6}{c|}{Runtime (s)} \\
\hline
Test Case & {SOC} & {SOC+LNC} & {CDF} & {CDF+LNC} & {NF} & {NF+LNC} & {SOC} & {SOC+LNC} & {CDF} & {CDF+LNC} & {NF} & {NF+LNC} \\
\hline
case14-ieee-sad        &    6.05 &    3.45 &    6.05 &    3.45 &   14.70 &    8.99 &     4 ms &     4 ms &     2 ms &     3 ms &     1 ms &     2 ms \\
case24-ieee-rts-sad    &    5.88 &    2.38 &    5.88 &    2.38 &   17.76 &   17.76 &     7 ms &     7 ms &     7 ms &     7 ms &     1 ms &     2 ms \\
case30-ieee-sad        &    5.66 &    5.33 &    5.66 &    5.33 &    6.89 &    5.96 &     6 ms &     7 ms &     7 ms &     7 ms &     1 ms &     2 ms \\
case30-as-sad          &   19.26 &    5.34 &   19.26 &    5.34 &   46.00 &   28.67 &     7 ms &     8 ms &     9 ms &     7 ms &     1 ms &     2 ms \\
case39-epri-sad        &    6.03 &    0.97 &    6.03 &    0.97 &   14.71 &   14.71 &  0.01    &  0.01    &  0.01    &     9 ms &     1 ms &     2 ms \\
case57-ieee-sad        &    4.88 &    2.44 &    4.88 &    2.44 &   40.74 &   37.51 &  0.01    &  0.02    &  0.01    &  0.02    &     3 ms &     4 ms \\
case60-c-sad           &   11.36 &    1.91 &   11.36 &    1.91 &   84.32 &   84.32 &  0.02    &  0.02    &  0.02    &  0.02    &     2 ms &     3 ms \\
case73-ieee-rts-sad    &    6.54 &    3.04 &    6.54 &    3.04 &   17.66 &   17.66 &  0.02    &  0.02    &  0.02    &  0.02    &     3 ms &     6 ms \\
case118-ieee-sad       &   19.90 &    8.42 &   19.90 &    8.42 &   45.59 &   45.59 &  0.04    &  0.04    &  0.04    &  0.04    &     5 ms &     7 ms \\
case300-ieee-sad       &    3.33 &    2.30 &    3.33 &    2.29 &   38.13 &   38.11 &  0.10    &  0.11    &  0.12    &  0.13    &  0.01    &  0.03    \\
case793-goc-sad        &   29.40 &   19.78 &   29.40 &   19.77 &   72.54 &   71.55 &  0.25    &  0.27    &  0.16    &  0.20    &  0.02    &  0.05    \\
case2312-goc-sad       &   19.41 &   17.69 &   19.41 &   17.69 &   53.91 &   53.83 &  1.06    &  1.24    &  0.61    &  1.02    &  0.08    &  0.20    \\
case3022-goc-sad       &   25.63 &   19.19 &   25.65 &   19.20 &  127.49 &  127.48 &  1.47    &  1.75    &  0.80    &  0.88    &  0.08    &  0.18    \\
\hline
case14-ieee-api        &   16.74 &    1.64 &   16.74 &    1.64 &   23.86 &   23.86 &     4 ms &     5 ms &     3 ms &     3 ms &     1 ms &     1 ms \\
case24-ieee-rts-api    &    6.58 &    1.17 &    6.58 &    1.17 &   37.88 &   37.59 &     7 ms &     9 ms &     8 ms &     8 ms &     1 ms &     1 ms \\
case30-ieee-api        &   14.51 &    1.11 &   14.51 &    1.11 &   20.26 &   18.47 &     7 ms &     7 ms &     7 ms &     8 ms &     1 ms &     2 ms \\
case30-as-api          &    9.87 &    0.40 &    9.87 &    0.40 &   50.22 &   50.22 &     7 ms &     9 ms &     8 ms &     8 ms &     1 ms &     2 ms \\
case39-epri-api        &    1.47 &    0.36 &    1.47 &    0.36 &    5.19 &    5.19 &  0.01    &  0.01    &     9 ms &     9 ms &     2 ms &     1 ms \\
case57-ieee-api        &   24.24 &    6.79 &   24.24 &    6.79 &   96.62 &   94.31 &  0.01    &  0.02    &  0.02    &  0.02    &     2 ms &     3 ms \\
case60-c-api           &    3.91 &    1.14 &    3.91 &    1.14 &   15.18 &   15.18 &  0.02    &  0.02    &  0.02    &  0.02    &     3 ms &     3 ms \\
case73-ieee-rts-api    &    6.97 &    1.92 &    6.97 &    1.92 &   42.43 &   41.74 &  0.02    &  0.03    &  0.03    &  0.03    &     2 ms &     4 ms \\
case118-ieee-api       &   15.31 &    9.01 &   15.31 &    9.01 &   17.35 &   17.18 &  0.04    &  0.05    &  0.04    &  0.04    &     4 ms &  0.01    \\
case300-ieee-api       &    7.63 &    4.74 &    7.63 &    4.72 &   25.08 &   25.08 &  0.10    &  0.12    &  0.14    &  0.13    &     9 ms &  0.03    \\
case793-goc-api        &   24.70 &   18.41 &   24.70 &   18.41 &   40.72 &   40.72 &  0.24    &  0.28    &  0.16    &  0.19    &  0.02    &  0.03    \\
case2312-goc-api       &   17.01 &   16.56 &   17.01 &   16.56 &   19.16 &   19.16 &  0.89    &  1.17    &  0.59    &  0.61    &  0.06    &  0.10    \\
case3022-goc-api       &   23.00 &   18.41 &   23.00 &   18.41 &   36.66 &   36.66 &  1.47    &  1.69    &  0.76    &  0.86    &  0.09    &  0.18    \\
\hline
\end{tabular}
\end{adjustbox}
\label{table1}
\end{table*}

(1) The LNC strengthen the SOC and CDF relaxations in 100\% of the test cases, with
average and maximum differences in the optimality gaps of 6.2\% and 17.5\% respectively.
The NF relaxation is strengthened in 46.2\% of test cases, with average and maximum
differences in the optimality gaps of 1.3\% and 17.3\% respectively.

(2) The effect of the LNC on the optimality gaps is more pronounced when solving the
SOC and CDF relaxations. This is because the LNC is intended to strengthen the region
defined by Eq. (\ref{eq0}d), which is not present in the NF formulation.

(3) Table \ref{table2} emphasizes the runtime performance difference between the three
strengthened relaxations. Even though the SOC and CDF relaxations provide the same
relaxation quality, the strengthened CDF solves faster than the strengthened SOC.

(4) Coffrin et al. \cite{Coffrin2016} demonstrated that the NF relaxation is scalable
for large datasets due to its linearity. Here, the NF relaxation still shows good scalability,
even with the inclusion of the LNC, making it suitable for finding tighter optimality
gaps when the use of stronger relaxations is computationally prohibitive. It shows
appropriate scalability up to 78,484 buses, making it a good choice for obtaining fast
lower bounds in global solution algorithms for large networks.

{
We note that decreases in optimality gaps from 14\% to 1\%, as we observe with case30-ieee-api,
may have a significant impact on the runtime performance of iterative algorithms for robust
optimization such as that proposed by Molzahn and Roald \cite{Molzahn2018}. These authors
report that their method takes two to five iterations to converge, where each iteration
involves the solution to a convex relaxation of AC-OPF. If the LNCs can reduce the iteration
count of such an algorithm by one, this would correspond to a 20-50\% reduction in runtime.
}

\begin{table*}[h!]
\captionsetup{justification=centering} 
\caption{Performance comparison for a sample of the largest datasets in the PGLib-OPF library.}
\centering
\begin{adjustbox}{max width=\textwidth}
\begin{tabular}{|l|*{3}{S}|*{3}{S}|}
\hline
& \multicolumn{3}{c|}{\% Optimality Gap} & \multicolumn{3}{c|}{Runtime (s)} \\
\hline
Test Case & {SOC+LNC} & {CDF+LNC} & {NF+LNC} & {SOC+LNC} & {CDF+LNC} & {NF+LNC} \\
\hline
case4837-goc-sad       &  14.07 &  14.07 &  14.07 &    2 &    1 &    1 \\
case5658-epigrids-sad  &  26.21 &  26.20 &  35.60 &    4 &    3 &    1 \\
case9591-goc-sad       &   4.15 &   4.15 &   4.15 &    4 &    4 &    2 \\
case24464-goc-sad      &  15.42 &  15.42 &  15.48 &   37 &   14 &    4 \\
case30000-goc-sad      &  20.44 &  20.44 &  43.02 &   32 &   14 &    2 \\
case78484-epigrids-sad &  25.01 &  25.01 &  35.10 &  151 &   54 &   11 \\
\hline
case4837-goc-api       &  16.10 &  16.10 &  30.83 &    3 &    1 &  0.3 \\
case5658-epigrids-api  &  29.63 &  29.63 &  42.29 &    4 &    2 &  0.3 \\
case9591-goc-api       &  19.05 &  19.05 &  30.45 &   10 &    4 &    1 \\
case24464-goc-api      &  27.60 &  27.59 &  34.04 &   28 &   10 &    3 \\
case30000-goc-api      &  18.05 &  18.05 &  36.08 &   29 &   10 &    1 \\
case78484-epigrids-api &  30.91 &  30.91 &  41.37 &  135 &   43 &    9 \\
\hline
\end{tabular}
\end{adjustbox}
\label{table2}
\end{table*}

\section{Conclusion}
This letter demonstrates that the projection of lifted nonlinear cuts into the variable
space of the SOC, CDF and NF relaxations has the potential to produce tighter optimality
gaps with minimal additional runtime overheads. We showed the trade-off between relaxation
quality and solve time, concluding that even though the strengthened SOC and CDF formulations
are equivalent, the strengthened CDF is the better alternative for solving large datasets.
While the NF relaxation provides a weaker optimality gap than CDF, it could be a better
choice for computing fast lower bounds during branch and bound algorithms for datasets
with more than 78,484 buses.    

\section{Acknowledgements}
We acknowledge funding from the Los Alamos National Laboratory LDRD program through the Center
for Nonlinear Studies (CNLS).
LA-UR-24-23960.


 

\bibliographystyle{IEEEtran}

\bibliography{References.bib}











\end{document}